\documentclass[letterpaper, 10 pt, conference]{ieeeconf}  

\IEEEoverridecommandlockouts                              
\usepackage{graphicx}
\usepackage{fixltx2e}
\usepackage{stfloats}

\usepackage{epsfig} 
\usepackage{epstopdf}

\usepackage{amsmath} 
\usepackage{amssymb}  

\usepackage{cite}
\usepackage{url}
\usepackage{hyperref}
\usepackage{booktabs}

\newcommand{\matpower}{M{\sc atpower}~}

\title{\LARGE \bf
Solution of Optimal Power Flow Problems using Moment Relaxations Augmented with Objective Function Penalization
}

\author{Daniel Molzahn, C\'{e}dric Josz, Ian Hiskens, and Patrick Panciatici
\thanks{Support from the Dow Sustainability Fellowship, ARPA-E grant \mbox{DE-AR0000232}, and Los Alamos National Laboratory subcontract 270958.}
\thanks{D. Molzahn and I. Hiskens are with the Dept. of Electrical Engineering and Computer Science,
        University of Michigan, Ann Arbor, MI 48109, USA
        {\tt\small molzahn@umich.edu, hiskens@umich.edu}}%
\thanks{C. Josz and P. Panciatici are with the R\&D Dept. of RTE, Versailles, France
        {\tt\small cedric.josz@rte-france.com, patrick.panciatici@rte-france.com }}%
}

\begin{document}

\maketitle
\thispagestyle{empty}
\pagestyle{empty}

\begin{abstract}
The optimal power flow (OPF) problem minimizes the operating cost of an electric power system. Applications of convex relaxation techniques to the non-convex OPF problem have been of recent interest, including work using the Lasserre hierarchy of ``moment'' relaxations to globally solve many OPF problems. By preprocessing the network model to eliminate low-impedance lines, this paper demonstrates the capability of the moment relaxations to globally solve large OPF problems that minimize active power losses for portions of several European power systems. Large problems with more general objective functions have thus far been computationally intractable for current formulations of the moment relaxations. To overcome this limitation, this paper proposes the combination of an objective function penalization with the moment relaxations. This combination yields feasible points with objective function values that are close to the global optimum of several large OPF problems. Compared to an existing penalization method, the combination of penalization and the moment relaxations eliminates the need to specify one of the penalty parameters and solves a broader class of problems.
\end{abstract}

\section{Introduction}

Determining the most efficient operating point for a power system requires solution of the optimal power flow (OPF) problem. This problem optimizes power system operation subject to both network equality constraints and engineering limits. Non-linearity of the constraint equations generally makes the OPF problem non-convex and can lead to local optima~\cite{bukhsh_tps}. Many solution techniques have been proposed, including successive quadratic programs, Lagrangian relaxation, heuristic optimization, and interior point methods~\cite{opf_litreview1993IandII,ferc4}. Some of these techniques calculate at-least-locally optimal solutions to many large problems. However, while many local solution techniques often find global solutions~\cite{molzahn_lesieutre_demarco-global_optimality_condition}, they may fail to converge or converge to a local optimum~\cite{bukhsh_tps,ferc5}.

Recent research attention has focused on convex relaxations of the OPF problem, which provide lower bounds on the optimal objective value and can certify infeasibility. Further, a convex relaxation based on semidefinite programming (SDP) globally solves many OPF problems~\cite{lavaei_tps}.

The SDP relaxation of~\cite{lavaei_tps} has been generalized to a family of ``moment relaxations'' using the Lasserre hierarchy~\cite{lasserre_book} for polynomial optimization~\cite{pscc2014,patrick,ibm_paper}. The moment relaxations take the form of SDPs, and the first-order relaxation in this hierarchy is equivalent to the SDP relaxation of~\cite{lavaei_tps}. Increasing the relaxation order in this hierarchy enables global solution of a broader class of OPF problems.

The ability to globally solve a broader class of OPF problems has a computational cost; the moment relaxations quickly become intractable with increasing order. Fortunately, second- and third-order moment relaxations globally solve many small problems for which the first-order relaxation fails to yield the globally optimal decision variables.

However, increasing system size results in computational challenges even for low-order moment relaxations. The second-order relaxation is computationally intractable for OPF problems with more than ten buses. Exploiting network sparsity enables solution of the first-order relaxation for systems with thousands of buses~\cite{jabr11,molzahn_holzer_lesieutre_demarco-large_scale_sdp_opf} and the second-order relaxation for systems with about forty buses~\cite{molzahn_hiskens-sparse_moment_opf,ibm_paper}. Recent work~\cite{molzahn_hiskens-sparse_moment_opf} solves larger problems (up to 300 buses) by both exploiting sparsity and only applying the computationally intensive higher-order moment relaxations to specific buses in the network. Other recent work improves computational tractability using a second-order cone programming relaxation of the higher-order moment constraints~\cite{powertech2015}.

Solving larger problems using moment relaxations is often limited by numerical convergence issues rather than computational performance. We present a preprocessing method that improves numerical convergence by removing low-impedance lines from the network model. Similar methods are commonly employed (e.g., PSS/E~\cite{PSSEManual}), but more extensive modifications are needed for adequate convergence due to the limited numerical capabilities of current SDP solvers.

After this preprocessing, the moment relaxations globally solve several large OPF problems which minimize active power losses for European power systems. Directly using the moment relaxations to globally solve more general large OPF problems with objective functions that minimize generation cost has been computationally intractable thus far.

To solve these OPF problems, we form moment relaxations using a penalized objective function. Previous literature~\cite{lavaei_mesh,lavaei_allerton2014} augments the SDP relaxation~\cite{lavaei_tps} with penalization terms for the total reactive power generation and the apparent power loss of certain lines. For many problems, this penalization finds feasible points with objective function values that are very close to the lower bounds obtained from the SDP relaxation. Related work~\cite{molzahn_josz_hiskens_panciatici-Laplacian_Objective} uses a Laplacian-based objective function with a constraint on generation cost to find feasible points are very near the global optima. This paper analyzes the physical and convexity properties of the reactive power penalization.

There are several disadvantages of the penalization method of~\cite{lavaei_allerton2014}. This penalization often requires choosing multiple parameters. (See~\cite{molzahn_josz_hiskens_panciatici-Laplacian_Objective} for a related approach that does not require choosing penalty parameters.) Also, there are OPF problems that are globally solved by the moment relaxations, but no known penalty parameters yield a feasible solution.

We propose a ``moment+penalization'' approach that augments the moment relaxations with a reactive power penalty. Typical penalized OPF problems only require higher-order moment constraints at a few buses. Thus, for a variety of large test cases, augmenting the moment relaxation with the proposed single-parameter penalization achieves feasible solutions that are at least very near the global optima (within at least 1\% for a variety of example problems). The moment+penalization approach enables solution of a broader class of problems than either method individually.

This paper is organized as follows. Section~\ref{l:opf_formulation} introduces the OPF problem. Section~\ref{l:moment} reviews the moment relaxations. Section~\ref{l:preprocess} describes the low-impedance line preprocessing. Section~\ref{l:penalization} discusses the existing penalization and the proposed moment+penalization approaches. Section~\ref{l:results} demonstrates the moment+penalization approach using several large test cases, and Section~\ref{l:conclusion} concludes the paper.

\section{Optimal Power Flow Problem}
\label{l:opf_formulation}

We first present an OPF formulation in terms of rectangular voltage coordinates, active and reactive power injections, and apparent power line-flow limits. Consider an $n$-bus system, where $\mathcal{N} = \left\lbrace 1, \ldots, n \right\rbrace$ is the set of buses, $\mathcal{G}$ is the set of generator buses, and $\mathcal{L}$ is the set of lines. Let $P_{Dk} + j Q_{Dk}$ represent the active and reactive load demand and $V_k = V_{dk} + j V_{qk}$ the voltage phasors at each bus~$k \in \mathcal{N}$. Superscripts ``max'' and ``min'' denote specified upper and lower limits. Buses without generators have maximum and minimum generation set to zero. Denote the network admittance matrix as $\mathbf{Y} = \mathbf{G} + j \mathbf{B}$.

Define a convex quadratic cost of active power generation:

\vspace{-8pt}
\begin{equation}\label{objfunction}\small
f_{Ck}\left(V_d,V_q\right) = c_{k2} \left(f_{Pk}\left(V_d,V_q\right)\right)^2 + c_{k1} f_{Pk}\left(V_d,V_q\right) + c_{k0}
\end{equation}
\vspace{-8pt}

The power flow equations describe the network physics:

\vspace{-10pt}
\begin{subequations}
\small
\label{opf_balance}
\begin{align}\nonumber
P_{Gk} = & f_{Pk}\left(V_d,V_q\right) = V_{dk} \sum_{i=1}^n \left( \mathbf{G}_{ik} V_{di} - \mathbf{B}_{ik} V_{qi} \right) &  &  \\[-5pt] 
\label{opf_Pbalance}  & + V_{qk} \sum_{i=1}^n \left( \mathbf{B}_{ik}V_{di} + \mathbf{G}_{ik}V_{qi} \right) + P_{Dk} \\ \nonumber
Q_{Gk} = & f_{Qk}\left(V_d,V_q \right) = V_{dk} \sum_{i=1}^n \left( -\mathbf{B}_{ik}V_{di} - \mathbf{G}_{ik} V_{qi}\right) \\[-5pt]
\label{opf_Qbalance} & + V_{qk} \sum_{i=1}^n \left( \mathbf{G}_{ik} V_{di} - \mathbf{B}_{ik} V_{qi}\right) + Q_{Dk}
\end{align}
\end{subequations}

Define a function for squared voltage magnitude:

\vspace{-8pt}
\begin{equation} \label{opf_Vsq}\small
\left(V_{k}\right)^2 = f_{Vk}\left(V_d, V_q\right) = V_{dk}^2 + V_{qk}^2
\end{equation}

We use a line model with an ideal transformer that has a specified turns ratio $\tau_{lm} e^{j\theta_{lm}} \colon 1$ in series with a $\Pi$ circuit with series impedance $R_{lm} + j X_{lm}$ (equivalent to an admittance of $g_{lm} + j b_{lm} = \frac{1}{R_{lm} + j X_{lm}}$) and shunt admittance $j b_{sh,lm}$. The line-flow equations are

\vspace{-8pt}
\begin{subequations}
\small
\begin{align}
\nonumber  & P_{lm} = f_{Plm}\left(V_d,V_q\right) = \left( V_{dl}^2 + V_{ql}^2\right) \frac{g_{lm}}{\tau_{lm}^2} \\ \nonumber
& \quad +
\left(V_{dl}V_{dm} + V_{ql}V_{qm}\right) \left(b_{lm}\sin\left(\theta_{lm} \right) - g_{lm}\cos\left(\theta_{lm} \right) \right) / \tau_{lm} \\
\label{Plm}& \quad +
\left(V_{dl}V_{qm} - V_{ql}V_{dm}\right) \left(g_{lm}\sin\left(\theta_{lm}\right) + b_{lm}\cos\left(\theta_{lm}\right)\right) / \tau_{lm} \\
\nonumber & P_{ml} = f_{Pml}\left(V_d,V_q\right) = \left(V_{dm}^2 + V_{qm}^2 \right)g_{lm} \\ \nonumber & \quad -
\left(V_{dl}V_{dm} + V_{ql}V_{qm} \right)\left(g_{lm}\cos\left(\theta_{lm}\right) + b_{lm}\sin\left(\theta_{lm}\right) \right) / \tau_{lm} \\ \label{Pml}
& \quad + 
\left(V_{dl}V_{qm} - V_{ql}V_{dm} \right) \left(g_{lm}\sin\left(\theta_{lm}\right) - b_{lm}\cos\left(\theta_{lm}\right) \right) / \tau_{lm}\\
\nonumber & Q_{lm} = f_{Qlm}\left(V_d,V_q\right) = -\left( V_{dl}^2 + V_{ql}^2\right) \left(b_{lm} + \frac{b_{sh,lm}}{2}\right) / \tau_{lm}^2 \\ \nonumber 
& \quad +
\left(V_{dl}V_{dm} + V_{ql}V_{qm}\right) \left(b_{lm}\cos\left(\theta_{lm} \right) + g_{lm}\sin\left(\theta_{lm} \right) \right) / \tau_{lm} \\ \label{Qlm} & \quad +
\left(V_{dl}V_{qm} - V_{ql}V_{dm}\right) \left(g_{lm}\cos\left(\theta_{lm}\right) - b_{lm}\sin\left(\theta_{lm}\right)\right) / \tau_{lm} \\
\nonumber & Q_{ml} = f_{Qml}\left(V_d,V_q\right) = -\left( V_{dm}^2 + V_{qm}^2\right) \left(b_{lm} + \frac{b_{sh,lm}}{2}\right) \\ \nonumber & \quad +
\left(V_{dl}V_{dm} + V_{ql}V_{qm}\right) \left(b_{lm}\cos\left(\theta_{lm} \right) - g_{lm}\sin\left(\theta_{lm} \right) \right)  / \tau_{lm} \\ \label{Qml} & \quad + 
\left(-V_{dl}V_{qm} + V_{ql}V_{dm}\right) \left(g_{lm}\cos\left(\theta_{lm}\right) + b_{lm}\sin\left(\theta_{lm}\right)\right) / \tau_{lm} \\
\label{Slm} & \left(S_{lm}\right)^2 = f_{Slm}\left(V_d,V_q\right) = \left(f_{Plm}\left(V_d,V_q\right)\right)^2 + \left(f_{Qlm}\left(V_d,V_q\right)\right)^2 \\
\label{Sml} & \left(S_{ml}\right)^2 = f_{Sml}\left(V_d,V_q\right) = \left(f_{Pml}\left(V_d,V_q\right)\right)^2 + \left(f_{Qml}\left(V_d,V_q\right)\right)^2
\end{align}
\end{subequations}

The classical OPF problem is then

\vspace{-8pt}
\begin{subequations}
\label{opf}
\small
\begin{align}
\label{opf_obj} & \min_{V_d,V_q} \sum_{k \in \mathcal{G}} f_{Ck}\left(V_d,V_q\right) \qquad\qquad\quad \mathrm{subject\; to} \hspace{-60pt} & \\
\label{opf_P} &  \quad P_{Gk}^{\mathrm{min}} \leq f_{Pk}\left(V_d,V_q\right) \leq P_{Gk}^{\mathrm{max}} & \forall k \in \mathcal{N} \\
\label{opf_Q} &  \quad Q_{Gk}^{\mathrm{min}} \leq f_{Qk}\left(V_d,V_q\right) \leq Q_{Gk}^{\mathrm{max}} &  \forall k \in \mathcal{N} \\
\label{opf_V} &  \quad (V_{k}^{\mathrm{min}})^2 \leq f_{Vk}\left(V_d,V_q\right) \leq (V_{k}^{\mathrm{max}})^2 &  \forall k \in \mathcal{N}  \\
\label{opf_Slm} & \quad f_{Slm}\left(V_d,V_q\right) \leq \left(S_{lm}^{\mathrm{max}}\right)^2 &  \forall \left(l,m\right) \in \mathcal{L} \\ 
\label{opf_Sml} & \quad f_{Sml}\left(V_d,V_q\right) \leq \left(S_{lm}^{\mathrm{max}}\right)^2 &  \forall \left(l,m\right) \in \mathcal{L} \\ 
\label{opf_Vref} & \quad V_{q1} = 0
\end{align}
\end{subequations}
 
Constraint~\eqref{opf_Vref} sets the reference bus angle to zero.

\section{Moment Relaxations}
\label{l:moment}

Since all constraints and the objective function are polynomials in the voltage components $V_d$ and $V_q$, the OPF problem \eqref{opf} is a polynomial optimization problem that can be solved with tools from algebraic geometry. We first review the ``moment'' relaxations from the Lasserre hierarchy~\cite{lasserre_book} for polynomial optimization problems and then summarize a method for exploiting network sparsity.

\subsection{Review of Moment Relaxations}
\label{l:momentreview}

Polynomial optimization problems are a special case of ``generalized moment problems''~\cite{lasserre_book}. Global solutions to generalized moment problems can be approximated using moment relaxations that are formulated as SDPs. For polynomial optimization problems with bounded variables, such as OPF problems, the approximation approaches the global solution(s) as the relaxation order increases~\cite{lasserre_book}. The first-order relaxation in the Lasserre hierarchy is equivalent to the SDP relaxation of~\cite{lavaei_tps}; higher-order moment relaxations generalize the SDP relaxation of~\cite{lavaei_tps}.

We begin with several definitions. Let the vector $\hat{x} = \begin{bmatrix} V_{d1} & V_{d2} & \ldots & V_{qn} \end{bmatrix}^\intercal$ contain all first-order monomials of the decision variables in~\eqref{opf}. Given a vector \mbox{$\alpha \in \mathbb{N}^{2n}$} representing monomial exponents, the expression $\hat{x}^\alpha = V_{d1}^{\alpha_1}V_{d2}^{\alpha_2}\cdots V_{qn}^{\alpha_{2n}}$ defines the monomial associated with $\hat{x}$ and $\alpha$. A polynomial $g\left(\hat{x}\right)$ is

\vspace{-5pt}
\begin{equation}\label{gpoly}
g\left(\hat{x}\right) \triangleq \sum_{\alpha \in \mathbb{N}^{2n}} g_{\alpha} \hat{x}^{\alpha}
\end{equation}
\vspace{-8pt}

\noindent where $g_{\alpha}$ is the scalar coefficient corresponding to $\hat{x}^{\alpha}$.

Define a linear functional $L_y\left\lbrace g\right\rbrace$ which replaces the monomials $\hat{x}^\alpha$ in a polynomial $g\left(\hat{x}\right)$ with scalars $y_{\alpha}$:

\vspace{-5pt}
\begin{equation}\label{L}
L_y\left\lbrace g\right\rbrace \triangleq \sum_{\alpha \in \mathbb{N}^{2n}} g_{\alpha} y_{\alpha}
\end{equation}
\vspace{-8pt}

\noindent Apply $L_y\left\lbrace g \right\rbrace$ to each element of a matrix argument.

Consider the vector $\hat{x} = \begin{bmatrix}V_{d1} & V_{d2} & V_{q2} \end{bmatrix}^\intercal$ containing the voltage components of a two-bus system, where the angle reference~\eqref{opf_Vref} is used to eliminate $V_{q1}$, and the polynomial $g\left(\hat{x}\right) = -\left(0.95\right)^2 + f_{V2}\left(V_d,V_q \right) = -\left(0.95\right)^2 + V_{d2}^2 + V_{q2}^2$. (Constraining $g\left(\hat{x}\right) \geq 0$ forces the voltage magnitude at bus~2 to be greater than or equal to 0.95~per unit.) Then $L_y\left\lbrace g\right\rbrace = -\left(0.95\right)^2y_{000} + y_{020} + y_{002}$. Thus, $L_y\left\lbrace g \right\rbrace$ converts a polynomial $g\left(\hat{x}\right)$ to a linear function of $y$.

Define the vector $x_\gamma$ composed of all monomials of the voltage components up to the specified relaxation order $\gamma$:

\vspace{-5pt}
\begin{align} \nonumber
x_\gamma \triangleq & \left[ \begin{array}{ccccccc} 1 & V_{d1} & \ldots & V_{qn} & V_{d1}^2 & V_{d1}V_{d2} & \ldots \end{array} \right. \\ \label{x_d}
& \qquad \left.\begin{array}{cccccc} \ldots & V_{qn}^2 & V_{d1}^3 & V_{d1}^2 V_{d2} & \ldots & V_{qn}^\gamma \end{array}\right]^\intercal
\end{align}
\vspace{-8pt}

The moment relaxations constrain \emph{moment} and \emph{localizing} matrices. The symmetric moment matrix $\mathbf{M}_\gamma\left\lbrace y\right\rbrace$ has entries $y_\alpha$ corresponding to all monomials $\hat{x}^\alpha$ up to order $2\gamma$:

\vspace{-5pt}
\begin{equation}\label{momentmat}
\mathbf{M}_\gamma \left\lbrace y \right\rbrace \triangleq L_y\left\lbrace x_\gamma^{\vphantom{\intercal}} x_\gamma^\intercal\right\rbrace
\end{equation}
\vspace{-8pt}

Symmetric localizing matrices are defined for each constraint of~\eqref{opf}. The localizing matrices consist of linear combinations of the moment matrix entries $y$. Each polynomial constraint of the form $f\left(\hat{x}\right) - a \geq 0$ in~\eqref{opf} (e.g., $f_{V2}\left(\hat{x}\right) - V_2^{\min} \geq 0$) corresponds to the localizing matrix

\vspace{-5pt}
\begin{equation}\label{localizing} \small
\mathbf{M}_{\gamma-\beta}\left\lbrace \left(f\left(\hat{x}\right) - a\right) y \right\rbrace \triangleq L_y\left\lbrace\left(f\left(\hat{x}\right) - a\right) x_{\gamma-\beta}^{\vphantom{\intercal}} x_{\gamma-\beta}^\intercal \right\rbrace
\end{equation}
\vspace{-8pt}

\noindent where the polynomial $f$ has degree $2\beta$. Example moment and localizing matrices for the second-order relaxation of a two-bus system are presented in~\eqref{2busMomentMat} and~\eqref{2busLocalizingMat}, respectively.

The order-$\gamma$ moment relaxation of~\eqref{opf} is

\vspace{-5pt}
\begin{subequations}\small
\label{msdp_opf}
\begin{align}
\label{msdp_obj}& \min_{y} L_y\left\lbrace \sum_{k \in \mathcal{G}} f_{Ck} \right\rbrace \qquad \mathrm{subject\; to} \hspace{-150pt} &  \\
\label{msdp_Pmin} & \quad \mathbf{M}_{\gamma-1}\left\lbrace \left(f_{Pk} - P_k^{\min}\right) y \right\rbrace \succeq 0 & \forall k\in\mathcal{N}\\
\label{msdp_Pmax} & \quad \mathbf{M}_{\gamma-1}\left\lbrace \left(P_k^{\max} - f_{Pk} \vphantom{P_k^{\min}}\right) y \right\rbrace \succeq 0 & \forall k\in\mathcal{N}\\
\label{msdp_Qmin} & \quad \mathbf{M}_{\gamma-1}\left\lbrace \left(f_{Qk} - Q_k^{\min}\right) y \right\rbrace \succeq 0 & \forall k\in\mathcal{N}\\
\label{msdp_Qmax} & \quad \mathbf{M}_{\gamma-1}\left\lbrace \left(Q_k^{\max} - f_{Qk}  \vphantom{P_k^{\min}}\right) y \right\rbrace \succeq 0 & \forall k\in\mathcal{N}\\
\label{msdp_Vmin} & \quad \mathbf{M}_{\gamma-1}\left\lbrace \left(f_{Vk} - V_k^{\min}\right) y \right\rbrace \succeq 0 & \forall k\in\mathcal{N}
\end{align}\begin{align}
\label{msdp_Vmax} & \quad \mathbf{M}_{\gamma-1}\left\lbrace \left(V_k^{\max} - f_{Vk}  \vphantom{P_k^{\min}}\right) y \right\rbrace \succeq 0 & \forall k\in\mathcal{N} \\
\label{msdp_Slm} & \quad \mathbf{M}_{\gamma-2}\left\lbrace \left(S_{lm}^{\max} - f_{Slm} \vphantom{P_k^{\min}}\right) y \right\rbrace \succeq 0 & \forall \left(l,m\right)\in\mathcal{L} \\
\label{msdp_Sml} & \quad \mathbf{M}_{\gamma-2}\left\lbrace \left(S_{lm}^{\max} - f_{Sml} \vphantom{P_k^{\min}}\right) y \right\rbrace \succeq 0 & \forall \left(l,m\right)\in\mathcal{L} \\
\label{msdp_Msdp} & \quad \mathbf{M}_\gamma \left\lbrace y\right\rbrace \succeq 0 & \\
\label{msdp_y0} & \quad y_{00\ldots 0} = 1 & \\
\label{msdp_Vref} & \quad y_{0\ldots00\eta0\ldots0} = 0 & \eta = 1,\ldots,2\gamma
\end{align}
\end{subequations}

\begin{figure*}[t]
\small
\setcounter{equation}{12}
\begin{equation}\label{2busMomentX} 
x_2 = \left[\begin{array}{cccccccccc} 1 & V_{d1} & V_{d2} & V_{q2} & V_{d1}^2 & V_{d1}V_{d2} & V_{d1}V_{q2} & V_{d2}^2 & V_{d2}V_{q2} & V_{q2}^2\end{array}\right]^\intercal \qquad\; \left[\mathrm{Note\!:}\; \eqref{opf_Vref}\; \mathrm{is\; used\; to\; remove}\; V_{q1} \right]
\end{equation}
\vspace{-6pt}
\begin{equation}\label{2busMomentMat}\small
\mathbf{M}_2 \left\lbrace y \right\rbrace = L_y\left\lbrace x_2^{\vphantom{\intercal}} x_2^\intercal\right\rbrace = \left[\begin{array}{c|ccc|cccccc} 
y_{000} & y_{100} & y_{010} & y_{001} & y_{200} & y_{110} & y_{101} & y_{020} & y_{011} & y_{002} \\\hline
y_{100} & y_{200} & y_{110} & y_{101} & y_{300} & y_{210} & y_{201} & y_{120} & y_{111} & y_{102} \\
y_{010} & y_{110} & y_{020} & y_{011} & y_{210} & y_{120} & y_{111} & y_{030} & y_{021} & y_{012} \\
y_{001} & y_{101} & y_{011} & y_{002} & y_{201} & y_{111} & y_{102} & y_{021} & y_{012} & y_{003} \\ \hline
y_{200} & y_{300} & y_{210} & y_{201} & y_{400} & y_{310} & y_{301} & y_{220} & y_{211} & y_{202} \\ 
y_{110} & y_{210} & y_{120} & y_{111} & y_{310} & y_{220} & y_{211} & y_{130} & y_{121} & y_{112} \\
y_{101} & y_{201} & y_{111} & y_{102} & y_{301} & y_{211} & y_{202} & y_{121} & y_{112} & y_{103} \\
y_{020} & y_{120} & y_{030} & y_{021} & y_{220} & y_{130} & y_{121} & y_{040} & y_{031} & y_{022} \\
y_{011} & y_{111} & y_{021} & y_{012} & y_{211} & y_{121} & y_{112} & y_{031} & y_{022} & y_{013} \\ 
y_{002} & y_{102} & y_{012} & y_{003} & y_{202} & y_{112} & y_{103} & y_{022} & y_{013} & y_{004} 
\end{array}\right]
\end{equation}
\vspace{-13pt}
\begin{align}\label{2busLocalizingMat}  \small\nonumber
& \mathbf{M}_{1}\left\lbrace\left(f_{V2} - \left(0.95\right)^2\right) y \right\rbrace = \\
& \quad \left[\begin{array}{c|ccc}
y_{020} + y_{002} - \left(0.95\right)^2y_{000} & y_{120} + y_{102} - \left(0.95\right)^2y_{100} & y_{030} + y_{012} - \left(0.95\right)^2y_{010} & y_{021} + y_{003} - \left(0.95\right)^2y_{001}  \\ \hline
y_{120} + y_{102} - \left(0.95\right)^2y_{100} & y_{220} + y_{202} - \left(0.95\right)^2y_{200} & y_{130} + y_{112} - \left(0.95\right)^2y_{110} & y_{121} + y_{103} - \left(0.95\right)^2y_{101} \\
y_{030} + y_{012} - \left(0.95\right)^2y_{010} & y_{130} + y_{112} - \left(0.95\right)^2y_{110} & y_{040} + y_{022} - \left(0.95\right)^2y_{020} & y_{031} + y_{013} - \left(0.95\right)^2y_{011} \\
y_{021} + y_{003} - \left(0.95\right)^2y_{001} & y_{121} + y_{103} - \left(0.95\right)^2y_{101} & y_{031} + y_{013} - \left(0.95\right)^2y_{011} & y_{022} + y_{004} - \left(0.95\right)^2y_{002}  \\
\end{array}\right]
\end{align}
\setcounter{equation}{11}
\vspace*{0pt}
\hrule
\vspace*{-6pt}
\end{figure*}

\noindent where $\succeq 0$ indicates that the corresponding matrix is positive semidefinite. The constraint~\eqref{msdp_y0} enforces $x^{0} = 1$. The constraint~\eqref{msdp_Vref} corresponds to the angle reference~\eqref{opf_Vref}; the $\eta$ in \eqref{msdp_Vref} is in the index $n+1$, which corresponds to the variable $V_{q1}$. The angle reference can alternatively be used to eliminate all terms corresponding to $V_{q1}$ to reduce the size of the SDP, as in~\eqref{2busMomentX}--\eqref{2busLocalizingMat}.

The objective function and apparent power line flow constraints are quartic polynomials in $V_d$ and $V_q$. For $\gamma = 1$, use a Schur complement formulation of these polynomials~\cite{molzahn_hiskens-sparse_moment_opf}.

The \mbox{order-$\gamma$} moment relaxation yields a single global solution upon satisfaction of the rank condition 

\vspace{-6pt}
\begin{equation}\label{rankcondition}
\mathrm{rank}\left(\mathbf{M}_{\gamma}\left\lbrace y\right\rbrace\right) = 1
\end{equation}
\vspace{-10pt}

\noindent The global solution $x^\ast$ to the OPF problem~\eqref{opf} is then determined by a spectral decomposition of the diagonal block of the moment matrix corresponding to the second-order terms. Specifically, let $\mu$ be a unit-length eigenvector corresponding to the non-zero eigenvalue $\lambda$ from the diagonal block of the moment matrix corresponding to the second-order monomials (i.e., $\left[\mathbf{M}_\gamma\left\lbrace y \right\rbrace\right]_{2:k,2:k}$, where $k=2n+1$ and subscripts indicate entries in MATLAB notation). Then the vector $V^\ast = \sqrt{\lambda} \left(\mu_{1:n} + j \mu_{\left(n+1\right):2n}\right)$ is the global solution.

A solution with $\mathrm{rank}\left(\mathbf{M}_{\gamma}\left\lbrace y\right\rbrace \right) > 1$ indicates that the order-$\gamma$ moment relaxation only yields a lower bound on the objective value. Increasing the relaxation order will improve the lower bound and may give a global solution.

\subsection{Exploiting Network Sparsity}
\label{l:sparsemoment}

The matrices in the moment relaxations quickly grow with both the relaxation order and the system size. For an $n$-bus system, the number of rows and columns in the \mbox{order-$\gamma$} relaxation's moment matrix is $\left(2n+\gamma\right)! / \left( \left(2n\right)! \gamma!\right)$. Solving second- and higher-order moment relaxations of problems with more than about ten buses requires exploiting network sparsity~\cite{ibm_paper,molzahn_hiskens-sparse_moment_opf}. Full details are excluded from this summary for brevity; see~\cite{molzahn_hiskens-sparse_moment_opf} for a complete description.

We use a matrix completion theorem~\cite{gron1984} to exploit network sparsity. A symmetric matrix $\mathbf{W}$ with partial information (i.e., not all entries of $\mathbf{W}$ have known values) can be completed to a positive semidefinite matrix (i.e., the unknown entries of $\mathbf{W}$ can be chosen such that $\mathbf{W} \succeq 0$) if and only if certain submatrices of $\mathbf{W}$ are positive semidefinite. Consider a \emph{chordal extension} of the power system network graph.\footnote{A chordal extension contains all links in the network as well as additional links such that every cycle of length four or more nodes has an edge connecting two non-adjacent nodes in the cycle. A chordal extension can be calculated using a Cholseky factorization of the network Laplacian matrix using an approximate minimize-degree permutation to maintain sparsity.} The matrix $\mathbf{W}$ is positive semidefinite if and only if the submatrices of $\mathbf{W}$ corresponding to all maximal cliques (i.e., the largest completely connected subgraphs) of the chordal extension are positive semidefinite.

The matrix completion theorem enables decomposition of a positive semidefinite constraint for a single large matrix to constraints on many smaller matrices. This eliminates many terms which do not appear in the constraint equations of the OPF problem~\eqref{opf}. See~\cite{waki2006,jabr11,molzahn_holzer_lesieutre_demarco-large_scale_sdp_opf,molzahn_hiskens-sparse_moment_opf} for details.

Exploiting network sparsity enables solution of the second-order relaxation for problems with up to about 40 buses. Recognizing that the first-order relaxation is sufficient for large regions of typical OPF problems, previous work~\cite{molzahn_hiskens-sparse_moment_opf} proposed an iterative algorithm to selectively apply the higher-order relaxation constraints to specific regions of the network~\cite{molzahn_hiskens-sparse_moment_opf}. This algorithm used a ``power injection mismatch'' heuristic to determine where to apply the higher-order moment constraints.

We next summarize this heuristic. A moment relaxation yields the global optimum if the moment matrix has rank one, in which case the algorithm terminates. If the rank of the moment matrix is greater than one, we calculate the closest rank-one matrix using an eigen decomposition. The power injection mismatch heuristic compares the power injections corresponding to the closest rank-one matrix with the power injections corresponding to the higher-rank moment matrix. Typical OPF problems result in small power injection mismatches at the majority of buses with a small fraction of buses having large mismatches. Each iteration of the algorithm in~\cite{molzahn_hiskens-sparse_moment_opf} tightens the relaxation by applying higher-order relaxation constraints for the two buses with greatest power injection mismatch.

\section{Preprocessing Low-Impedance Lines}
\label{l:preprocess}

By exploiting sparsity and selectively applying the higher-order constraints, the moment relaxations globally solve many OPF problems with up to 300 buses. Solution of larger problems with higher-order relaxations is typically limited by numerical convergence issues rather than computational concerns. This section describes a preprocessing method for improving numerical properties of the moment relaxations.

Low-impedance lines, which often represent ``jumpers'' between buses in the same physical location, cause numerical problems for many algorithms. Low line impedances result in a large range of values in the bus admittance matrix $\mathbf{Y}$, which causes numerical problems in the constraint equations.

To address these numerical problems, many software packages remove lines with impedances below a threshold. For instance, lines with impedance below a parameter \texttt{thrshz} are removed prior to applying other algorithms in PSS/E~\cite{PSSEManual}.

We use a slightly modified version of the low-impedance line removal procedure in PSS/E~\cite{PSSEManual}.\footnote{Lines with non-zero resistances are not considered to be ``low impedance'' by PSS/E. We consider both the resistance and the reactance.} Low-impedance lines are eliminated by grouping buses that are connected by lines with impedances below a specified threshold \texttt{thrshz}. Each group of buses is replaced by one bus that is connected to all lines terminating outside the group. Generators, loads, and shunts (including the shunt susceptanes of lines connecting buses within a group) are aggregated. The series parameters of lines connecting buses within a group are eliminated.

Removing low-impedance lines typically has a small impact on the solution. To recover an approximate solution to the original power system model, assign identical voltage phasors to all buses in each group and distribute flows on lines connecting buses within a group under the approximation that all power flows through the low-impedance lines.

A typical low-impedance line threshold \texttt{thrshz} is \mbox{$1\times 10^{-4}$}~per unit. However, the numerical capabilities of SDP solvers are not as mature as other optimization tools. Therefore, we require a larger $\mathtt{thrshz} = 1\times 10^{-3}$ per unit to obtain adequate convergence of the moment relaxations. This larger threshold typically introduces only small errors in the results, although non-negligible errors are possible.

\matpower solutions obtained for the Polish~\cite{matpower} and most \mbox{PEGASE} systems~\cite{pegase} were the the same before and after low-impedance line preprocessing to within 0.0095 per unit voltage magnitude and $0.67^\circ$ voltage angle difference across each line. Operating costs for all test problems were the same to within 0.04\%. The 2869-bus \mbox{PEGASE} system had larger differences: 0.0287 per unit voltage magnitude and $1.37^\circ$ angle difference. A power flow solution for the full network using the solution to the OPF problem after low-impedance line preprocessing yields smaller differences: 0.0059 per unit voltage magnitude and $1.17^\circ$ angle difference. Thus, the differences from the preprocessing for the 2869-bus \mbox{PEGASE} system can be largely attributed to the sensitivity of the OPF problem itself to small changes in the low-impedance line parameters. Preprocessing reduced the number of buses by between 21\% and 27\% for the PEGASE and between 9\% and 26\% for the Polish systems.

These results show the need for further study of the sensitivity of OPF problems to low-impedance line parameters as well as additional numerical improvements of the moment relaxations and SDP solvers to reduce \texttt{thrshz}. 


\section{Moment Relaxations and Penalization}
\label{l:penalization}

As will be shown in Section~\ref{l:results}, the moment relaxations globally solve many large OPF problems with active power loss minimization objectives after removing low-impedance lines as described in Section~\ref{l:preprocess}. Directly applying the moment relaxations to many large OPF problems with more general cost functions has so far been computationally intractable. This section describes the non-convexity associated with more general cost functions and proposes a method to obtain feasible solutions that are at least near the global optimum, if not, in fact, globally optimal for many problems.

Specifically, we propose augmenting the moment relaxations with a penalization in the objective function. Previous literature~\cite{lavaei_mesh,lavaei_allerton2014} adds terms to the first-order moment relaxation that penalize the total reactive power injection and the apparent power line loss (i.e., $\sqrt{\left(f_{Plm} + f_{Pml}\right)^2 + \left(f_{Qlm} + f_{Qml}\right)^2}$) for ``problematic'' lines identified by a heuristic. This penalization often finds feasible points that are at least nearly globally optimal.

However, the penalization in~\cite{lavaei_allerton2014} requires choosing two penalty parameters and fails to yield a feasible solution to some problems. This section describes progress in addressing these limitations by augmenting the moment relaxations with a reactive power penalization. The proposed ``moment+penalization'' approach only requires a single penalty parameter and finds feasible points that are at least nearly globally optimal for a broader class of OPF problems. This section also analyzes the convexity properties and provides a physical intuition for reactive power penalization.

\subsection{Penalization of Reactive Power Generation}
\label{l:penalty_explanation}

The penalization method proposed in~\cite{lavaei_allerton2014} perturbs the objective function~\eqref{objfunction} to include terms that minimize the total reactive power loss and the apparent power loss on specific lines determined by a heuristic method. These terms enter the objective function with two scalar parameter coefficients. Obtaining a feasible point near the global solution requires appropriate choice of these parameters.

For typical operating conditions, reactive power is strongly associated with voltage magnitude. Penalizing reactive power injections tends to reduce voltage magnitudes, which also tends to increase active power losses since a larger current flow, with higher associated loss, is required to deliver a given quantity of power at a lower voltage magnitude.

For many problems for which the first-order moment relaxation fails to yield the global optimum, the relaxation ``artificially'' increases the voltage magnitudes to reduce active power losses. This results in voltage magnitudes and power injections that are feasible for the relaxation~\eqref{msdp_opf} but infeasible for the OPF problem~\eqref{opf}. 

By choosing a reactive power penalty parameter that balances these competing tendencies (increasing voltage magnitudes to reduce active power losses vs. decreasing voltage magnitudes to reduce the penalty), the penalized relaxation finds a feasible solution to many OPF problems. Since losses typically account for a small percentage of active power generation and active and reactive power are typically loosely coupled, the reactive power penalization often results in a feasible point that is near the global optimum.

We next study the convexity properties of the cost function and the reactive power penalization. The cost function~\eqref{objfunction} is convex in terms of active power generation but not necessarily in terms of the voltage components $V_d$ and $V_q$ due to the non-convexity of~\eqref{opf_Pbalance}.\footnote{The cost function of the \emph{moment relaxation} \eqref{msdp_obj} is always convex. This section studies the convexity of the penalized objective function for the \emph{original} non-convex OPF problem~\eqref{opf_obj}.} Thus, the objective function is a potential source of non-convexity which may result in the relaxation's failure to globally solve the OPF problem.

Consider the eigenvalues of the symmetric matrices $\mathbf{C}$ and $\mathbf{D}$, where, for the vector $\hat{x}$ containing the voltage components, $\hat{x}^\intercal \mathbf{C} \hat{x}$ is a linear cost of active power generation (i.e., $c_2 = c_0 = 0$ in~\eqref{objfunction}) and $\hat{x}^\intercal \mathbf{D} \hat{x}$ calculates the reactive power losses. For the 2383-bus Polish system~\cite{matpower}, which has linear generation costs, the most negative eigenvalue of $\mathbf{C}$ is $-8.53\times 10^7$. Thus, the objective function of the original OPF problem~\eqref{opf_obj} is non-convex in terms of the voltage components, which can cause the relaxation~\eqref{msdp_opf} to fail to yield the global optimum. Conversely, active power loss minimization is convex in terms of the voltage components due to the absence of negative resistances.

As indicated by the potential for negative eigenvalues of $\mathbf{D}$ (e.g., the matrix $\mathbf{D}$ for the 2383-bus Polish system has a pair of negative eigenvalues at $-0.0175$), penalizing reactive power losses is generally non-convex due to capacitive network elements (i.e., increasing voltage magnitudes may decrease the reactive power loss). See~\cite{molzahn_josz_hiskens_panciatici-Laplacian_Objective} for related work that uses a convex objective based on a Laplacian matrix.


Further work is needed to investigate the effects of reactive power penalization on OPF problems with more realistic generator models that explicitly consider the trade-off between active and reactive power outputs (i.e., generator ``D-curves''). A tighter coupling between active and reactive power generation may cause the reactive power penalization to yield solutions that are far from the global optimum.

The apparent power line loss penalty's effects are not as easy to interpret as the reactive power penalty. Ongoing work includes understanding the effects of the line loss penalty.

\subsection{Moment+Penalization Approach}

Although the reactive power penalization often yields a near rank-one solution (i.e., rank-one matrices for most of the submatrices in the decomposition discussed in Section~\ref{l:sparsemoment}), this penalization alone is not sufficient to obtain a feasible point for many problems. Reference~\cite{lavaei_allerton2014} penalizes the apparent power line loss associated with certain lines to address the few remaining non-rank-one ``problematic'' submatrices. However, this approach has several disadvantages.

First, penalizing apparent power line losses introduces another parameter.\footnote{Reference~\cite{lavaei_allerton2014} uses the same penalization parameter for each ``problematic'' line. Generally, each line could have a different penalty parameter.} Introducing parameters is problematic, especially when lacking an intuition for appropriate values.

Second, the combination of reactive power and line loss penalization may not yield a feasible solution to some problems. For instance, the OPF problems case9mod and case39mod1 from~\cite{bukhsh_tps} are globally solved with low-order moment relaxations, but there is no known penalization of reactive power and/or apparent power line loss that yields a feasible solution for these problems. Also, the penalization approach is not guaranteed to yield a feasible solution that is close to the global optimum.

Unlike the penalization approach, the moment relaxation approach does not require the choice of penalty parameters, globally solves a broader class of OPF problems, and is guaranteed to yield the global optimum when the rank condition~\eqref{rankcondition} is satisfied. However, direct application of the moment relaxations to large problems has so far been limited to active power loss minimization objective functions. We conjecture that the non-convexity associated with more general cost functions requires higher-order moment constraints at too many buses for computational tractability.

To apply the moment relaxations to large OPF problems with active power generation cost objective functions, we augment the moment relaxations with a reactive power penalty. Specifically, we apply the sparsity-exploiting moment relaxations described in Section~\ref{l:moment} to the OPF problem~\eqref{opf} where the objective function~\eqref{opf_obj} is replaced by 

\vspace{-5pt}
\begin{align}\label{Qpenobj}
\sum_{k \in \mathcal{G}} \left(f_{Ck}\left(V_d,V_q\right) + \epsilon_b f_{Qk}\left(V_d,V_q\right)\right)
\end{align}
\vspace{-5pt}

\noindent where $\epsilon_b$ is the scalar reactive power penalization parameter. That is, rather than apply an apparent power loss penalization to the objective function, we apply higher-order moment constraints to specific buses~\cite{molzahn_hiskens-sparse_moment_opf}. As will be demonstrated in Section~\ref{l:results}, higher-order moment constraints are only needed at a few buses in typical OPF problems after augmenting the objective function with a reactive power penalization term.

Similar to the existing penalization, when the rank condition~\eqref{rankcondition} is satisfied, the proposed ``moment+penalization'' approach yields the global solution to the \emph{modified} OPF problem~\eqref{Qpenobj}, but not necessarily to the original OPF problem~\eqref{opf}. However, since the penalization does not change the constraint equations, the solution to the moment+penalization approach is \emph{feasible} for the original OPF problem~\eqref{Qpenobj}. The first-order moment relaxation without penalization (i.e., $\epsilon_b = 0$) gives a lower bound on the globally optimal objective value for the original OPF problem~\eqref{opf}. This lower bound provides an optimality metric for the feasible solution obtained from the moment+penalization approach. As will be shown in Section~\ref{l:results}, the feasible solutions for a variety of problems are within at least 1\% of the global optimum.

The moment+penalization approach inherits a mix of the advantages and disadvantages of the moment relaxation and penalization methods. First, the moment+penalization approach requires selection of a single scalar parameter (one more than needed for the moment relaxations, but one less than generally needed for the penalization in~\cite{lavaei_allerton2014}). This parameter must be large enough to result in a near rank-one solution (i.e., the solution should have many rank-one submatrices in the decomposition described in Section~\ref{l:sparsemoment}), but small enough to avoid large changes to the OPF problem.

Second, the penalization eliminates the moment relaxations' guarantees: the moment+penalization approach may yield a feasible solution that is far from the global optimum or not give any solution. However, the moment+penalization approach finds global or near-global solutions to a broader class of small OPF problems than penalization approach of~\cite{lavaei_allerton2014} (e.g., case9mod and case39mod1 with $\epsilon_b = 0$, and case39mod3 with $\epsilon_b = \$0.10/\mathrm{MVAr}$~\cite{bukhsh_tps}). This suggests that the moment+penalization approach inherits the ability of the moment relaxations to solve a broad class of OPF problems.

Finally, the penalization in the moment+penalization approach enables calculation of feasible solutions that are at least nearly globally optimal for a variety of large OPF problems with objective functions that minimize active power generation cost rather than just active power losses.

Note that it is not straightforward to compare the computational costs of the moment+penalization approach and the penalization approach in~\cite{lavaei_allerton2014}. A single solution of a penalized first-order moment relaxation, as in~\cite{lavaei_allerton2014}, is faster than a relaxation with higher-order moment constraints. Thus, if one knows appropriate penalty parameters, the method in~\cite{lavaei_allerton2014} is faster. Although a relatively wide range of penalty parameters tends to work well for typical OPF problems, there are problems for which no known penalty parameters yield feasible solutions. For these problems, the moment+penalization approach has a clear advantage.

The moment+penalization approach has the advantage of systematically tightening the relaxation rather than requiring the choice of penalty parameters. However, the higher-order constraints can significantly increase solver times. Thus, there is a potential trade-off between finding appropriate penalization parameters for the approach in~\cite{lavaei_allerton2014} and increased solver time from the moment+penalization approach. The speed of the moment+penalization approach may be improved using the mixed SDP/SOCP relaxation from~\cite{powertech2015}.

\section{Results}
\label{l:results}

This section first globally solves several large, active-power-loss minimizing OPF problems using moment relaxations without penalization (\mbox{$\epsilon_b = 0$}).  Next, this section applies the moment+penalization approach to find feasible points that are at least nearly globally optimal for several test cases which minimize active power generation cost. Unless otherwise stated, the preprocessing method from Section~\ref{l:preprocess} with \texttt{thrshz} set to $1\times 10^{-3}$ per unit is applied to all examples. No example enforces a minimum line resistance.

The results are generated using the iterative algorithm from~\cite{molzahn_hiskens-sparse_moment_opf} which selectively applies the higher-order moment relaxation constraints as summarized in Section~\ref{l:sparsemoment}. The algorithm terminates when all power injection mismatches are less than 1~MVA.

The implementation uses \mbox{MATLAB 2013a}, YALMIP 2015.06.26~\cite{yalmip}, and Mosek 7.1.0.28, and was solved using a computer with a quad-core 2.70~GHz processor and 16~GB of RAM. The test cases are the Polish system models in \matpower\cite{matpower} and several \mbox{PEGASE} systems~\cite{pegase} representing portions of the European power system.


\subsection{Active Power Loss Minimization Results}

Table~\ref{t:lossresults} shows the results of applying the moment relaxations to several large OPF problems that minimize active power losses (i.e., the cost coefficients in~\eqref{objfunction} are $c_2 = c_0 = 0$, $c_1 = \$1/\textrm{MWh}$). The solutions to the preprocessed problems are guaranteed to be globally optimal since there is no penalization. The columns in Table~\ref{t:lossresults} list the case name, the number of iterations of the algorithm from~\cite{molzahn_hiskens-sparse_moment_opf}, the maximum power injection mismatch, the globally optimal objective value, and the solver time summed over all iterations.\footnote{PEGASE-1354 and PEGASE-2869 use a \texttt{thrshz} of \mbox{$3\times 10^{-3}$}~per unit. All other systems use $1\times 10^{-3}$ per unit.} The abbreviation ``PL'' stands for ``Poland''. Table~\ref{t:lossresults} excludes several cases (the 89-bus \mbox{PEGASE} system and the Polish 2736sp, 2737sop, and 2746wp systems) which only require the first-order relaxation and thus do not illustrate the capabilities of the higher-order relaxations.

\begin{table}[thb]
\caption{Active Power Loss Minimization Results}
\label{t:lossresults}
\footnotesize
\centering
\begin{tabular}{|c|c|c|c|c|}
\hline 
\textbf{Case} & \textbf{Num.}  & \textbf{Global Obj.} & \textbf{Max $\mathrm{S}^{\mathrm{mis}}$}& \textbf{Solver} \\
\textbf{Name} & \textbf{Iter.} & \textbf{Val. (\$/hr)} & \textbf{(MVA)} & \textbf{\!\!Time (sec)\!\!} \\ \hline
PL-2383wp & 3 & \hphantom{1}$24990$ & 0.25 & \hphantom{1}583\\ 
PL-2746wop & 2 & \hphantom{1}$19210$ & 0.39 & 2662\\ 
PL-3012wp & 5 & \hphantom{1}$27642$ & 1.00 & \hphantom{1}319\\ 
PL-3120sp & 7 & \hphantom{1}$21512$ & 0.77 & \hphantom{1}387\\ 
PEGASE-1354 & 5 & \hphantom{1}$74043$ & 0.85 & \hphantom{1}407 \\ 
PEGASE-2869 & 6 & $133944$ & 0.63 & \hphantom{1}921 \\\hline 
\end{tabular}
\end{table}

Each iteration of the algorithm in~\cite{molzahn_hiskens-sparse_moment_opf} after the first adds second-order constraints at two buses. Thus, a small number of second-order buses (between 0.1\% and 0.7\% of the number of buses in the systems in Table~\ref{t:lossresults} after the low-impedance line preprocessing) are applied to all examples in Table~\ref{t:lossresults}. This results in computational tractability for the moment relaxations.

Note that PL-2746wop has a much greater solver time than the other systems even though it only has second-order constraints at two buses. This slow solution time is due to the fact that the two second-order buses are contained in submatrices corresponding to cliques with 10 and 11 buses. The second-order constraints for these large submatrices dominate the solver time. The mixed SDP/SOCP relaxation in~\cite{powertech2015} may be particularly useful beneficial for such cases.

Since the low-impedance line preprocessing has been applied to these systems, the solutions do not exactly match the original OPF problems. \matpower\cite{matpower} solutions of the original problems have objective values that are slightly larger than the values in Table~\ref{t:lossresults} due to losses associated with the line resistances removed by the preprocessing.

After the low-impedance line preprocessing, local solutions from \matpower match the solutions from the moment relaxations and are therefore, in fact, globally optimal. This is not the case for all OPF problems~\cite{bukhsh_tps,molzahn_hiskens-sparse_moment_opf}.

\subsection{Moment+Penalization for More General Cost Functions}

As discussed in Section~\ref{l:penalty_explanation}, minimization of active power generation cost often yields a non-convex objective function in terms of the voltage components. Despite this non-convexity, low-order moment relaxations typically yield global solutions to small problems, including problems without known penalty parameters for obtaining a feasible points (e.g., case9mod and case39mod1 from~\cite{bukhsh_tps}).

\begin{figure}[t]
\centering
\includegraphics[totalheight=0.18\textheight]{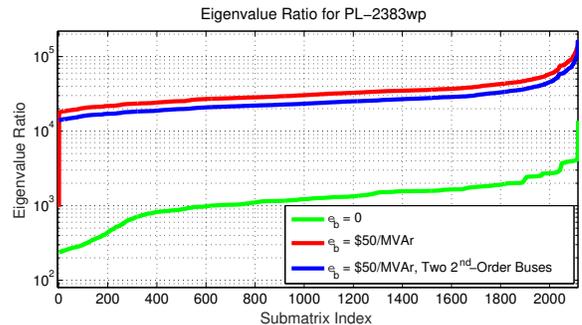}
\vspace{-20pt}
\caption{\hspace{-8pt} Eigenvalue ratio (largest/second-largest eigenvalue) for each submatrix for PL-2383wp. Large values ($>\!\!10^4$) indicate satisfaction of the rank condition~\eqref{rankcondition}. For $\epsilon_b = 0$ (green), most of the submatrices are not rank one. For $\epsilon_b = \$50/\mathrm{MVAr}$ (red), most submatrices satisfy the rank condition with the exception of those on the far left of the figure. Applying second-order moment constraints to two of the buses that are in these submatrices (blue) results in all submatrices satisfying the rank condition.}
\label{f:case2383wp_eigWratio}
\vspace{-13pt}
\end{figure}

However, the moment relaxations are thus far intractable for some large OPF problems with non-convex objective functions. A reactive power penalty often results in the first-order moment relaxation yielding a solution that is nearly globally optimal (i.e., most of the submatrices in the decomposition described in Section~\ref{l:sparsemoment}  satisfy the rank condition~\eqref{rankcondition}). Enforcing higher-order constraints at buses in the remaining submatrices yields a feasible solution to the OPF problem. This is illustrated in Fig.~\ref{f:case2383wp_eigWratio}, which shows the ratio between the largest and second-largest eigenvalues of the submatrices of the moment matrix, arranged in increasing order, for the 2383-bus Polish system. If the submatrices were all rank one, then this eigenvalue ratio would be infinite. Thus, large numeric values (i.e., greater than $1\times 10^4$) indicate satisfaction of the rank condition within numerical precision. Without the reactive power penalty, the rank condition is not satisfied for most submatrices. With the reactive power penalty, the rank condition is satisfied for many but not all submatrices. Enforcing higher-order moment constraints at two buses which are in the high-rank submatrices results in a feasible (rank-one) operating point for the OPF problem which is within at least 0.74\% of the global optimum.

To further illustrate the effectiveness of the moment+penalization approach, Table~\ref{t:costresults} shows the results of applying the moment+penalization approach to several large OPF problems with active power generation cost functions. The optimality gap column gives the percent difference between a lower bound on the optimal objective value from the first-order moment relaxation and the feasible solution obtained from the moment+penalization approach for the system after low-impedance line preprocessing.

\begin{table}[thb]
\caption{Generation Cost Minimization Results}
\label{t:costresults}
\footnotesize
\centering
\begin{tabular}{|c|c|c|c|c|c|}
\hline 
\textbf{Case} & \textbf{$\epsilon_b$} & \textbf{Num.}  & \textbf{Opt.} & \textbf{Max $\mathrm{S}^{\mathrm{mis}}$}& \textbf{Solver} \\
\textbf{Name} & \textbf{(\$/MVAr)} & \textbf{Iter.} & \textbf{Gap} & \textbf{(MVA)} & \textbf{Time (sec)} \\ \hline
PL-2383wp & 50 & 2 & 0.74\% & 0.13 & 152.2\\ 
PL-3012wp & 50 & 7 & 0.49\% & 0.20 & 1056.3\\ 
PL-3120sp & 100 & 6 & 0.92\% & 0.08 & 1164.4\\\hline
\end{tabular}
\end{table}

The penalized first-order relaxation requires 74.6, 88.9, and 97.0 seconds for PL-2383wp, PL-3012wp, and PL-3120sp, respectively. Attributing the rest of the solver time to the higher-order constraints implies that these constraints accounted for 3.1, 433.7, and 582.4 seconds beyond the time required to repeatedly solve the first-order relaxations.

The moment+penalization approach can yield feasible points that are at least nearly globally optimal for cases where both the penalization method of~\cite{lavaei_allerton2014} and low-order moment relaxations fail individually. For instance, the moment+penalization approach with a reactive power penalty of $\epsilon_b = \$0.10/\mathrm{MVAr}$ gives a feasible point within 0.28\% of the global optimum for case39mod3 from~\cite{bukhsh_tps}, but both second- and third-order moment relaxations and the penalization method in~\cite{lavaei_allerton2014} fail to yield global solutions.

\section{Conclusion}
\label{l:conclusion}

``Moment'' relaxations from the Lasserre hierarchy for polynomial optimization globally solve a broad class of OPF problems. By exploiting sparsity and selectively applying the computationally intensive higher-order moment relaxations, previous literature demonstrated the moment relaxations' capability to globally solve moderate-size OPF problems. This paper presented a preprocessing method that removes low-impedance lines to improve the numerical conditioning of the moment relaxations. After applying the preprocessing method, the moment relaxations globally solve a variety of OPF problems that minimize active power losses for systems with thousands of buses. A proposed ``moment+penalization'' method is capable of finding feasible points that are at least nearly globally optimal for large OPF problems with more general cost functions. This method has several advantages over previous penalization approaches, including requiring fewer parameter choices and solving a broader class of OPF problems. The results are demonstrated using large OPF problems representing European power systems. 



\bibliographystyle{IEEEtran}
\bibliography{IEEEabrv,cdc2015}

\begin{thebibliography}{10}
\providecommand{\url}[1]{#1}
\csname url@samestyle\endcsname
\providecommand{\newblock}{\relax}
\providecommand{\bibinfo}[2]{#2}
\providecommand{\BIBentrySTDinterwordspacing}{\spaceskip=0pt\relax}
\providecommand{\BIBentryALTinterwordstretchfactor}{4}
\providecommand{\BIBentryALTinterwordspacing}{\spaceskip=\fontdimen2\font plus
\BIBentryALTinterwordstretchfactor\fontdimen3\font minus
  \fontdimen4\font\relax}
\providecommand{\BIBforeignlanguage}[2]{{%
\expandafter\ifx\csname l@#1\endcsname\relax
\typeout{** WARNING: IEEEtran.bst: No hyphenation pattern has been}%
\typeout{** loaded for the language `#1'. Using the pattern for}%
\typeout{** the default language instead.}%
\else
\language=\csname l@#1\endcsname
\fi
#2}}
\providecommand{\BIBdecl}{\relax}
\BIBdecl

\bibitem{bukhsh_tps}
W.~Bukhsh, A.~Grothey, K.~McKinnon, and P.~Trodden, ``{Local Solutions of the
  Optimal Power Flow Problem},'' \emph{{IEEE} Trans. Power Syst.}, vol.~28,
  no.~4, pp. 4780--4788, 2013.

\bibitem{opf_litreview1993IandII}
J.~Momoh, R.~Adapa, and M.~El-Hawary, ``{A Review of Selected Optimal Power
  Flow Literature to 1993. Parts I and II},'' \emph{{IEEE} Trans. Power Syst.},
  vol.~14, no.~1, pp. 96--111, Feb. 1999.

\bibitem{ferc4}
A.~Castillo and R.~O'Neill, ``{Survey of Approaches to Solving the ACOPF (OPF
  Paper 4)},'' US Federal Energy Regulatory Commission, Tech. Rep., Mar. 2013.

\bibitem{molzahn_lesieutre_demarco-global_optimality_condition}
D.~Molzahn, B.~Lesieutre, and C.~DeMarco, ``{A Sufficient Condition for Global
  Optimality of Solutions to the Optimal Power Flow Problem},'' \emph{{IEEE}
  Trans. Power Syst.}, vol.~29, no.~2, pp. 978--979, Mar. 2014.

\bibitem{ferc5}
A.~Castillo and R.~O'Neill, ``{Computational Performance of Solution Techniques
  Applied to the ACOPF (OPF Paper 5)},'' US Federal Energy Regulatory
  Commission, Tech. Rep., Jan. 2013.

\bibitem{lavaei_tps}
J.~Lavaei and S.~Low, ``{Zero Duality Gap in Optimal Power Flow Problem},''
  \emph{{IEEE} Trans. Power Syst.}, vol.~27, no.~1, pp. 92--107, Feb. 2012.

\bibitem{lasserre_book}
J.-B. Lasserre, \emph{{Moments, Positive Polynomials and Their
  Applications}}.\hskip 1em plus 0.5em minus 0.4em\relax Imperial College
  Press, 2010, vol.~1.

\bibitem{pscc2014}
D.~Molzahn and I.~Hiskens, ``{Moment-Based Relaxation of the Optimal Power Flow
  Problem},'' \emph{18th Power Syst. Comput. Conf. (PSCC)}, 18-22 Aug. 2014.

\bibitem{patrick}
C.~Josz, J.~Maeght, P.~Panciatici, and J.~Gilbert, ``{Application of the
  Moment-SOS Approach to Global Optimization of the OPF Problem},''
  \emph{{IEEE} Trans. Power Syst.}, vol.~30, no.~1, pp. 463--470, Jan. 2015.

\bibitem{ibm_paper}
B.~Ghaddar, J.~Marecek, and M.~Mevissen, ``{Optimal Power Flow as a Polynomial
  Optimization Problem},'' \emph{\normalfont{To appear in} \textit{IEEE Trans.
  Power Syst.}}

\bibitem{jabr11}
R.~Jabr, ``{Exploiting Sparsity in SDP Relaxations of the OPF Problem},''
  \emph{{IEEE} Trans. Power Syst.}, vol.~27, no.~2, pp. 1138--1139, May 2012.

\bibitem{molzahn_holzer_lesieutre_demarco-large_scale_sdp_opf}
D.~Molzahn, J.~Holzer, B.~Lesieutre, and C.~DeMarco, ``{Implementation of a
  Large-Scale Optimal Power Flow Solver Based on Semidefinite Programming},''
  \emph{{IEEE} Trans. Power Syst.}, vol.~28, no.~4, pp. 3987--3998, 2013.

\bibitem{molzahn_hiskens-sparse_moment_opf}
D.~Molzahn and I.~Hiskens, ``{Sparsity-Exploiting Moment-Based Relaxations of
  the Optimal Power Flow Problem},'' \emph{IEEE Trans. Power Syst.}, vol.~30,
  no.~6, pp. 3168--3180, Nov. 2015.

\bibitem{powertech2015}
------, ``{Mixed SDP/SOCP Moment Relaxations of the Optimal Power Flow
  Problem},'' in \emph{IEEE Eindhoven PowerTech}, Jun. 2015.

\bibitem{PSSEManual}
{Siemens PTI}, ``{Volume II: Program Application Guide},'' \emph{Power System
  Simulation for Engineering (PSS/E)}, vol. 31.0, {December} 2007.

\bibitem{lavaei_mesh}
R.~Madani, S.~Sojoudi, and J.~Lavaei, ``{Convex Relaxation for Optimal Power
  Flow Problem: Mesh Networks},'' \emph{{IEEE} Trans. Power Syst.}, vol.~30,
  no.~1, pp. 199--211, Jan. 2015.

\bibitem{lavaei_allerton2014}
R.~Madani, M.~Ashraphijuo, and J.~Lavaei, ``{Promises of Conic Relaxation for
  Contingency-Constrained Optimal Power Flow Problem},'' in \emph{{52nd Annu.
  Allerton Conf. Commun., Control, and Comput.}}, Sept. 2014, pp. 1064--1071.

\bibitem{molzahn_josz_hiskens_panciatici-Laplacian_Objective}
D.~Molzahn, C.~Josz, I.~Hiskens, and P.~Panciatici, ``{A Laplacian-Based
  Approach for Finding Near Globally Optimal Solutions to OPF Problems},''
  \emph{Submitted. \rm{Preprint available:
  \url{http://arxiv.org/abs/1507.07212}}}.

\bibitem{gron1984}
R.~Gron, C.~Johnson, E.~S\'{a}, and H.~Wolkowicz, ``{Positive Definite
  Completions of Partial Hermitian Matrices},'' \emph{{Linear Algebra Appl.}},
  vol.~58, pp. 109--124, 1984.

\bibitem{waki2006}
H.~Waki, S.~Kim, M.~Kojima, and M.~Muramatsu, ``{Sums of Squares and
  Semidefinite Program Relaxations for Polynomial Optimization Problems with
  Structured Sparsity},'' \emph{SIAM J. Optimiz.}, vol.~17, no.~1, pp.
  218--242, 2006.

\bibitem{matpower}
R.~Zimmerman, C.~Murillo-S{\'a}nchez, and R.~Thomas, ``{MATPOWER: Steady-State
  Operations, Planning, and Analysis Tools for Power Systems Research and
  Education},'' \emph{{IEEE} Trans. Power Syst.}, no.~99, pp. 1--8, 2011.

\bibitem{pegase}
S.~Fliscounakis, P.~Panciatici, F.~Capitanescu, and L.~Wehenkel, ``{Contingency
  Ranking with Respect to Overloads in Very Large Power Systems Taking into
  Account Uncertainty, Preventive and Corrective Actions},'' \emph{{IEEE}
  Trans. Power Syst.}, vol.~28, no.~4, pp. 4909--4917, 2013, {(System models
  available in M{\sc atpower})}.

\bibitem{yalmip}
J.~L\"{o}fberg, ``{YALMIP: A Toolbox for Modeling and Optimization in
  MATLAB},'' in \emph{{IEEE Int. Symp. Compu. Aided Control Syst. Des.}}, 2004,
  pp. 284--289.

\end{thebibliography}

\end{document}